\def\vers{July 3, 2006}
\magnification=1200
\hsize=6.5truein
\vsize=8.9truein
\font\bigfont=cmr10 at 14pt
\font\sfont=cmr8
\font\sbfont=cmbx8
\font\sifont=cmti8
\def\ssbull{\raise.2ex\hbox{${\scriptscriptstyle\bullet}$}}
\def\mopls{\hbox{$\bigoplus$}}
\def\mcup{\hbox{$\bigcup$}}
\def\msum{\hbox{$\sum$}}
\def\wdg{\hbox{$\wedge$}}
\def\bC{{\bf C}}
\def\bN{{\bf N}}
\def\bQ{{\bf Q}}
\def\bZ{{\bf Z}}
\def\boR{{\bf R}}
\def\cA{{\cal A}}
\def\cB{{\cal B}}
\def\cC{{\cal C}}
\def\cD{{\cal D}}
\def\cG{{\cal G}}
\def\cH{{\cal H}}
\def\cK{{\cal K}}
\def\cO{{\cal O}}
\def\ocA{\bar{\cal A}}
\def\ocK{\bar{\cal K}}
\def\tcA{\widetilde{\cal A}}
\def\tD{\widetilde{D}}
\def\tf{\widetilde{f}}
\def\th{\widetilde{h}}
\def\tH{\widetilde{H}}
\def\tX{\widetilde{X}}
\def\tbet{\widetilde{\beta}}
\def\tgam{\widetilde{\gamma}}
\def\tcG{\widetilde{\cal G}}
\def\tcH{\widetilde{\cal H}}
\def\hatopl{\widehat{\hbox{$\bigoplus$}}}
\def\Im{\hbox{\rm Im}}
\def\Ker{\hbox{\rm Ker}}
\def\Coker{\hbox{\rm Coker}}
\def\can{\hbox{\rm can}}
\def\DR{\hbox{\rm DR}}
\def\Gr{\hbox{\rm Gr}}
\def\Sing{\hbox{\rm Sing}}
\def\rd{\partial}
\def\lbb{\{\!\{}
\def\rbb{\}\!\}}
\def\simto{\buildrel\sim\over\longrightarrow}
\def\text{\hbox}

\hbox{}
\vskip 1cm

\centerline{\bigfont Brieskorn Modules and Gauss-Manin Systems}

\centerline{\bigfont for Non-isolated Hypersurface Singularities}

\bigskip
\centerline{Daniel Barlet$^{\dag}$ and Morihiko Saito$^{\dag\dag}$}
\footnote{}{\sfont \hskip -5pt
$^{\dag}$Universit\'e Henri Poincar\'e (Nancy I) et
Institut Universitaire de France}
\footnote{}{\sfont
Institut E. Cartan UHP/CNRS/INRIA, UMR 7502, B.P. 239}
\footnote{}{\sfont 54506 Vandoeuvre-les-Nancy Cedex, FRANCE}
\footnote{}{\sfont e-mail: barlet@iecn.u-nancy.fr}
\footnote{}{\sfont \hskip -10pt
$^{\dag\dag}$RIMS Kyoto University, Kyoto 606-8502 JAPAN}
\footnote{}{\sfont e-mail: msaito@kurims.kyoto-u.ac.jp}

\bigskip\bigskip
{\narrower\noindent
{\sbfont Abstract.} {\sfont
We study the Brieskorn modules associated to a germ of holomorphic
function with non-isolated singularities, and show that the Brieskorn
module has naturally a structure of a module over the ring of
microdifferential operators of nonpositive degree,
and that the kernel of the morphism to the Gauss-Manin system
coincides with the torsion part for the action of $t$ and also with
that for the action of the inverse of the Gauss-Manin connection.
This torsion part is not finitely generated in general,
and we give a sufficient condition for the finiteness.
We also prove a Thom-Sebastiani type theorem for the sheaf of
Brieskorn modules in the case one of two functions has an
isolated singularity.}\par}

\bigskip\bigskip
\centerline{\bf Introduction}
\footnote{}{{\sifont Date}{\sfont : \vers}}
\footnote{}{{\sfont 2000}
{\sifont Mathematics Subject Classification.}
{\sfont 32S40}}

\bigskip\noindent
Let
$ f $ be a nonconstant holomorphic function on a complex
manifold
$ X $, and
$ x \in X_{0} := f^{-1}(0) $.
The notion of {\it Brieskorn module} associated to a Milnor
fibration of
$ f $ around
$ x $ was introduced by Brieskorn [4]
in the case of {\it isolated} hypersurface singularities in order
to calculated the Milnor cohomology and the monodromy in an
algebraic way.
The generalization to the case of {\it non-isolated} hypersurface
singularities was done by H. Hamm [10] (unpublished),
see also [9].
He showed the coherence modulo torsion of the relative de Rham
cohomology sheaf associated to the Milnor fibration (but the
torsion part is not finitely generated in general).

To the Milnor fibration, one can also associate the
{\it Gauss-Manin system} (see [16] for the isolated
singularity case).
This is by definition the
direct image of the structure sheaf
$ \cO_{X} $ as a
$ \cD_{X} $-module under the Milnor fibration.
It has been known that this Gauss-Manin system is always a
coherent (or more precisely, holonomic)
$ \cD $-module even in the non-isolated hypersurface
singularity case according to M. Kashiwara.
Furthermore, we can easily see that its stalk at the origin
defines a constructible {\it sheaf} on
$ X_{0} $ when
$ x \in X_{0} $ varies.
We have a natural morphism of the sheaf of the
$ i $-th Brieskorn modules
$ \cH_{f}^{i} $ to the sheaf of the
$ i $-th Gauss-Manin systems
$ \cG_{f}^{i} $ on
$ X_{0} $, and these are closely related to the vanishing
cohomology sheaf of degree
$ i-1 $ on
$ X_{0} $.
Using this morphism, we can get more precise information on the
Brieskorn modules, because the Gauss-Manin systems are easier
to describe.

Let
$ \bC\lbb\rd_{t}^{-1}\rbb $ be the ring of microdifferential
operators
$ \sum_{j\ge 0}a_{j}\rd_{t}^{-j} $ with constant coefficients
and nonpositive degree, satisfying the convergence condition
$ \sum_{j\ge 0}|a_{j}|r^{j}/j! < \infty $ for some
$ r > 0 $, see [11], [16].
Then both
$ \cH_{f}^{i} $ and
$ \cG_{f}^{i} $ have naturally a structure of a sheaf of
$ \bC\{t\} $-modules and of
$ \bC\lbb\rd_{t}^{-1}\rbb $-modules on
$ X_{0} $.
The latter structure for
$ \cH_{f}^{1} $ uses the integration without constant term,
which is similar to a canonical action of
$ \rd_{t}^{-1} $ on
$ \bC\{t\} $ via the indefinite integral from the origin,
see the proof of Proposition (2.4).
Let
$ V^{\alpha} $ denote the
$ V $-filtration of M. Kashiwara [12] and B. Malgrange [15] on
$ \cG_{f}^{i} $ indexed by
$ \bQ $.
For
$ \alpha > -1 $, it can be described by using the Deligne
extension [5] in this case.

In this paper we show the following (see (2.3--5)):

\medskip\noindent
{\bf Theorem 1.}
{\it The
$ t $-torsion part of
$ \cH_{f}^{i} $ coincide with the
$ \rd_{t}^{-1} $-torsion part, and also with the kernel of
$ \cH_{f}^{i} \to \cG_{f}^{i} $.
Furthermore these are annihilated by
$ t^{p} $ and
$ \rd_{t}^{-p} $ for
$ p \gg 0 $ locally on
$ X_{0} $.
The stalks of the image of the morphism
$ \cH_{f}^{i} \to \cG_{f}^{i} $ are finite free modules over
$ \bC\{t\} $ and over
$ \bC\lbb\rd_{t}^{-1}\rbb $.
The image is contained in
$ V^{>-1}\cG_{f}^{i} $ and its rank coincides with that
of the Milnor cohomology of degree
$ i-1 $ at each
point of
$ X_{0} $.
}

\medskip
This is a refinement of Hamm's result [10] mentioned above
(see also [9]),
and is closely related to positivity (or negativity) theorems
in [11], [14].
In general, the torsion part is not finitely generated,
see Remark (3.5).
However, using the theory of Kiehl-Verdier [13],
we can show its finiteness if there is a stratification
such that
$ f $ is locally trivial along its strata, see Theorem (3.3).
This was shown in [17].

In the case
$ \dim\Sing\,f = 1 $ and
$ f $ is locally trivial along
$ \Sing\,f \setminus \{0\} $, it is possible to prove this finiteness
using an elementary method, see [1].
There is a sufficient condition for the torsion-freeness of the
Brieskorn module, and we have a formula for the dimension of the Milnor
cohomology under this condition together with the above conditions on
$ \Sing\,f $, see [1].
In the case
$ \dim X = 2 $ these conditions are always satisfied so that
the formula is valid, see [1].
Furthermore, they are stable by replacing
$ f $ with a function
$ f + g $ on
$ X\times Y $ if
$ g $ is a holomorphic function with an isolated singularity on
$ Y $, see [1].
In the last situation we have also the Thom-Sebastiani type theorem
([19], [20]) for the sheaf of the reduced Brieskorn modules
(corresponding to reduced cohomology),
see Theorem (4.2) and also [1] for a special case.

Part of this work was done during a visit of the first author to
Kyoto, and he thanks RIMS for its hospitality.

\bigskip\bigskip
\centerline{\bf 1. Gauss-Manin systems}

\bigskip\noindent
{\bf 1.1.~Sheaf of Gauss-Manin systems.}
Let
$ f $ be a nonconstant holomorphic function on a complex
manifold
$ X $.
We may assume that the image of
$ f $ is contained in an open disk
$ S $ by restricting
$ X $ if necessary.
Let
$ t $ be the coordinate of
$ S $,
and put
$ \rd_{t} = \rd /\rd t $.
Let
$ i_{f} : X \to X\times S $ be the graph embedding by
$ f $.
Let
$ \cB_{f} $ be the direct image of
$ \cO_{X} $ by
$ i_{f} $ as a left
$ \cD_{X} $-module, see e.g. [11], [15], [16].
Then
$$
\cB_{f} = \cO_{X}\otimes_{\bC}\bC[\rd_{t}],
$$
and the action of a vector field
$ \xi $ on
$ X $ is given by
$$
\xi (\omega \otimes \rd_{t}^{i}) = \xi \omega \otimes
\rd_{t}^{i} - (\xi f)\omega \otimes \rd_{t}^{i+1}\quad
\text{for}\,\,\omega \in \Omega_{X}^{i}.
\leqno(1.1.1)
$$

Let
$$
\cK_{f}^{\ssbull} = \DR_{X\times S/S}(\cB_{f})[-\dim X].
$$
It is a complex whose
$ i $-th term is
$ \Omega_{X}^{i} \otimes_{\bC}
\bC[\rd_{t}] $ and whose differential is given by
$$
d(\omega \otimes \rd_{t}^{i}) = d\omega \otimes
\rd_{t}^{i} - df\wdg \omega \otimes \rd_{t}^{i+1}\quad
\text{for}\,\,\omega \in \Omega_{X}^{i}.
\leqno(1.1.2)
$$
Let
$ X_{0} = f^{-1}(0) $, and define
$$
\cG_{f}^{i} = \cH^{i}\cK_{f}^{\ssbull}|_{X_{0}},
$$
where
$ \cH^{i} $ denotes the cohomology sheaf of a sheaf complex.
We call
$ \cG_{f}^{i} $ the {\it sheaf of Gauss-Manin systems}
associated to
$ f $.
It has naturally a structure of
$ \cD_{S,0} $-modules, where
$ \cD_{S} $ is the sheaf of holomorphic linear differential
operators on
$ S $,
i.e.
$ \cD_{S,0} = \bC\{t\}\langle \rd_{t}\rangle $.

For
$ i = 1 $, we see that
$ \omega_{0} := df\otimes 1 \in \cH_{f}^{1} $ is annihilated by
$ \rd_{t} $ (i.e.
$ df\otimes \rd_{t} = 0 $ in
$ \cH_{f}^{1}) $, considering the image of
$ 1\otimes 1 $ by
$ d $.
So
$ \cD_{S,0}\omega_{0} = \bC\{t\}\omega_{0} $.
This is a free
$ \bC\{t\} $-module of rank
$ 1 $ using
$ [\rd_{t},t^{i}] = it^{i-1} $ inductively
(or by the involutivity of the characteristic variety.)
We define the {\it sheaf of reduced Gauss-Manin systems}
$ \tcG_{f}^{i} $ by
$$
\tcG_{f}^{i} = \cG_{f}^{i}\,\,\,
\hbox{\rm if}\,\,\, i \ne 1,\,\,\,\hbox{\rm and}\,\,\,
\cG_{f}^{1}/\bC\{t\}\omega_{0} \,\,\,
\hbox{\rm if}\,\,\, i = 1.
\leqno(1.1.3)
$$

\medskip\noindent
{\bf 1.2.~Nearby and vanishing cycles.}
We will denote by
$ \psi_{f}\bC_{X} $ and
$ \varphi_{f}\bC_{X} $ the nearby and vanishing cycle
sheaves, see [7].
Their cohomology sheaves
$ \cH^{i}\psi_{f}\bC_{X} $ and
$ \cH^{i}\varphi_{f}\bC_{X} $ are constructible sheaves of
$ \bC $-vector spaces.
They have naturally the action of the monodromy
$ T $.
By definition we have a distinguished triangle
$$
\bC_{X_{0}} \to \psi_{f}\bC_{X} \buildrel{\rm can}\over
\longrightarrow
\varphi_{f}\bC_{X} \buildrel{+1}\over\longrightarrow
\leqno(1.2.1)
$$

Let
$ F_{x} $ denote the Milnor fiber of
$ f $ around
$ x \in X_{0} $.
Then there are isomorphisms compatible with the action of the
monodromy
$ T $
$$
(\cH^{i}\psi_{f}\bC_{X})_{x} = H^{i}(F_{x},\bC),\quad
(\cH^{i}\varphi_{f}\bC_{X})_{x} = \tH^{i}(F_{x},\bC),
\leqno(1.2.2)
$$
where
$ \tH^{i}(F_{x},\bC) $ is the reduced cohomology.

\medskip\noindent
{\bf 1.3.~Regular holonomic
$ \cD $-modules in one variable.}
Let
$ M $ be a regular holonomic
$ \cD_{S,0} $-module.
For
$ \alpha \in \bC $,
let
$$
M^{\alpha} = \mcup_{i>0} \Ker((t\rd_{t} - \alpha)^{i} :
M \to M).
$$

We have the infinite direct sum decomposition
$$
M = \hatopl_{\alpha} M^{\alpha},
\leqno(1.3.1)
$$
where
$ \hatopl $ means the completion by an appropriate topology
(similar to the completion
$ \bC\{t\} $ of
$ \bC[t] $ or
$ \bC\lbb\rd_{t}^{-1}\rbb $ of
$ \bC[\rd_{t}^{-1}] $).
For simplicity, we assume that the direct sum is indexed by
$ \alpha \in \bQ $,
i.e.
$ M $ has quasi-unipotent monodromy.
We define the nearby and vanishing cycles by
$$
\psi_{t}M = \mopls_{-1<\alpha \le 0} M^{\alpha},\quad
\varphi_{t}M = \mopls_{-1\le \alpha <0} M^{\alpha}.
$$
It has the action of
$ T $ defined by
$ \exp(-2\pi it\rd_{t}) $ on
$ M^{\alpha} $.
We have the decompositions
$ \psi_{t}M = \psi_{t,1}M \oplus \psi_{t,\ne 1}M $,
$ \varphi_{t}M = \varphi_{t,1}M \oplus \varphi_{t,\ne 1}
M $ such that
$$
\psi_{t,1}M = M^{0},\quad \varphi_{t,1}M = M^{-1},\quad
\psi_{t,\ne 1}M = \varphi_{t,\ne 1}M =
\mopls_{-1<\alpha <0} M^{\alpha}.
$$
We have a morphism
$$
\can : \psi_{t}M \to \varphi_{t}M,
\leqno(1.3.2)
$$
whose restriction to
$ \psi_{t,1}M $ and
$ \psi_{t,\ne 1}M $ is given respectively by
$ \rd_{t} $ and the identity.
Note that
$ M $ is determined by
$ \psi_{t}M, \varphi_{t}M $ together with the action of
$ t, \rd_{t} $.

The filtration
$ V $ of Kashiwara [12] and Malgrange [15] is given by
$$
V^{\alpha}M = \hatopl_{\beta \ge \alpha} M^{\beta},
\quad V^{>\alpha}M = \hatopl_{\beta >\alpha} M^{\beta},
\leqno(1.3.3)
$$
and the above definition of
$ \psi, \varphi $ is compatible with the usual one.

We say that a
$ \bC\{t\} $-submodule
$ L $ of
$ M $ is a {\it lattice} if it is finite over
$ \bC\{t\} $ and generates
$ M $ over
$ \cD_{S,0} $.
The first condition is equivalent to
$ L\subset V^{\alpha} $ for some
$ \alpha \in \bQ $, and the second implies that
$ L\supset V^{\beta} $ for some
$ \beta \in \bQ $.
The equivalence holds for the last if the above morphism
$ \can : \psi_{t}M\to\varphi_{t}M $ is surjective;
for example, if the action of
$ \rd_{t} $ is bijective.
In the last case,
$ V^{\alpha}M $ and
$ M $ are finite free modules over
$ \bC\lbb\rd_{t}^{-1}\rbb $ and
$ \bC\lbb\rd_{t}^{-1}\rbb[\rd_{t}] $ respectively.
(This can be reduced to the fact that
$ \bC\{t\}t^{\alpha} \,(\alpha > 0) $ is a free
$ \lbb\rd_{t}^{-1}\rbb $-module of rank
$ 1 $ using a filtration.)
We have the notion of lattice with
$ \bC\{t\} $ replaced by
$ \bC\lbb\rd_{t}^{-1}\rbb $, and
$ \cD_{S,0} $ by
$ \bC\lbb\rd_{t}^{-1}\rbb[\rd_{t}] $ or by the ring of
microdifferential operators, see [11], [16].

Shrinking
$ S $ if necessary, we may assume that
$ M $ is the restriction of a regular holonomic
$ \cD_{S} $-module, which is also denoted by
$ M $,
and that the restriction of
$ M $ to
$ S\backslash \{0\} $ is coherent over
$ \cO_{S} $.
Note that
$ V^{\alpha}M $ for
$ \alpha > -1 $ is identified with the Deligne extension (see
[5]) such that the eigenvalues of the residue of
$ t\rd_{t} $ are contained in
$ [\alpha,\alpha +1) $.

The de Rham complex of
$ M $ is defined by
$$
\DR_{S}(M) = C(\rd_{t} : M \to M).
$$
It is well-known that we have canonical isomorphisms compatible
with the action of
$ T $
$$
\psi_{t}\DR_{S}(M)[-1] = \psi_{t}M,\quad \varphi_{t}
\DR_{S}(M)[-1] = \varphi_{t}M.
\leqno(1.3.4)
$$
(This is a special case of [12], [15].)
Indeed, the assertion for
$ \psi_{t} $ is quite classical (see [5] for a proof using
a modern language), and the assertion for
$ \varphi_{t} $ follows from it using a distinguished
triangle similar to (1.2.1).

Note that the argument in this subsection (1.3) can be applied to a
constructible sheaf of regular holonomic
$ \bC\{t\}\langle\rd_{t}\rangle $-modules on
$ X_{0} $, where
$ \bC\{t\}\langle\rd_{t}\rangle = \cD_{S,0} $.
In this case, each direct factor
$ M^{\alpha} $ of the decomposition (1.3.1) is a constructible
sheaf of
$ \bC $-vector spaces on
$ X_{0} $.

\medskip\noindent
{\bf 1.4.~Proposition.}
{\it The sheaves
$ \cG_{f}^{i}, \tcG_{f}^{i} $ are constructible sheaves of
regular holonomic
$ \cD_{S,0} $-modules, and there are canonical isomorphisms of
constructible sheaves of
$ \bC $-modules
$$
\psi_{t}\cG_{f}^{i+1} = \cH^{i}\psi_{f}\bC_{X},\quad
\varphi_{t}\cG_{f}^{i+1} = \cH^{i}\varphi_{f}\bC_{X},
\leqno(1.4.1)
$$
which are compatible with the action of
$ T $ and the morphism
$ \can $.
In particular,
$ \rd_{t} : \cG_{f}^{i} \to \cG_{f}^{i} $ is
bijective for
$ i \ne 1 $ and surjective for
$ i = 1 $ so that the kernel is the constant sheaf
$ \bC_{X_{0}} $.
On the reduced Gauss-Manin system
$ \tcG_{f}^{i} $, the action of
$ \rd_{t} $ is bijective for any
$ i $.
}

\medskip\noindent
{\it Proof.}
If
$ X_{0} $ is a divisor with normal crossings,
the assertion is well-known and easy to prove,
see also [18].
Let
$ \pi : \tX \to X $ be an embedded resolution of
$ X_{0} $,
i.e. the pull-back of
$ X_{0} $ is a divisor with simple normal crossings.
Then
$ \cO_{X} $ is a direct factor of the direct image
$ \pi_{+}\cO_{\tX} $ of
$ \cO_{\tX} $ as a
$ \cD $-module, and the complement is supported on
$ X_{0} $.
Indeed, we may assume that
$ \pi $ is an iteration of blowing-ups along
smooth centers by Hironaka, and the assertion is easily
verified for such morphisms.
(This is a special case of the decomposition theorem [2]
in the algebraic case.)

Let
$ \tf = f\pi $.
Then the above arguments imply that
$ \cK_{f}^{\ssbull} $ is a direct factor of
$ \boR\pi_{*}\cK_{\tf}^{\ssbull} $ and the
cohomology sheaves of its complement are
$ t $-torsion.
Furthermore,
$ \boR\pi_{*}\cK_{\tf}^{\ssbull}|_{X_{0}} $ has
constructible cohomology, using the spectral sequence
$$
E_{2}^{p,q} = R^{p}\pi_{*}(\cH^{q}\cK_{\tf}^{\ssbull}|_{\tX_{0}})
\Rightarrow \boR^{p+q}\pi_{*}\cK_{\tf}^{\ssbull}|_{X_{0}},
\leqno(1.4.2)
$$
which is compatible with the corresponding spectral sequences
for the nearby and vanishing cycles.
So
$ \cK_{f}^{\ssbull}|_{X_{0}} $ has constructible
cohomology sheaves of regular holonomic
$ \cD_{S,0} $-modules.
Thus we get (1.4.1) for
$ \psi_{t} $, and it remains to show the assertion
on the action of
$ \rd_{t} $ on
$ \cG_{f}^{i} $, because
$ \rd_{t} $ induces the morphism
$ \can : \psi_{t} \to \varphi_{t} $ in (1.3.2).

We have a short exact sequence
$$
0\to\cK_{f}^{\ssbull}\buildrel{\rd_{t}}\over\longrightarrow
\cK_{f}^{\ssbull}\to\Omega_{X}^{\ssbull} \to 0,
$$
and it induces a long exact sequence
$$
\to \cH^{i-1}\Omega_{X}^{\ssbull}|_{X_{0}} \to \cG_{f}^{i}
\buildrel{\rd_{t}}\over\longrightarrow \cG_{f}^{i}
\to \cH^{i}\Omega_{X}^{\ssbull}|_{X_{0}} \to,
$$
where
$ \cH^{i}\Omega_{X}^{\ssbull}|_{X_{0}} = \bC_{X_{0}} $ for
$ i = 0 $, and
$ \cH^{i}\Omega_{X}^{\ssbull}|_{X_{0}} = 0 $ otherwise.
So the assertion follows.

\medskip\noindent
{\bf 1.5.~Remark.}
Assume
$ f : X \to S $ is a Milnor representative at
$ x \in X_{0} $,
i.e.
$ X $ is the intersection of
$ f^{-1}(S) $ with an open
$ \varepsilon $-ball around
$ x $, and
$ S $ is an open disk of radius
$ \delta $ with
$ 0 < \delta \ll \varepsilon \ll 1 $.
Then it is known that the Gauss-Manin system
$ \boR^{i}f_{*}\cK_{f}^{\ssbull} $, which is the
direct image sheaf of
$ \cO_{X} $ by
$ f $ as a
$ \cD $-module, is a regular holonomic
$ \cD_{S} $-module and its stalk at the origin is
independent of
$ \varepsilon $.
This assertion follows from the theory of the direct
images of
$ \cD_{X} $-modules (see [11]), according to Kashiwara.
It can be verified also by using a resolution of
singularities
$ \pi : \tX\to X $ as in the proof of (1.4), because
we see that
$ \boR^{i}(f\pi)_{*}\cK_{\tf}^{\ssbull} $ is regular
holonomic for
$ 0 < \delta \ll \varepsilon \ll 1 $.

\medskip\noindent
{\bf 1.6.~Lemma.}
{\it For
$ i = 1 $, the short exact sequence
$$
0\to\bC\{t\}\omega_{0}\to\cG_{f}^{1}\to\tcG_{f}^{1}\to 0
\leqno(1.6.1)
$$
splits canonically, and we have in the notation of {\rm (1.3.1)}
}
$$
(\tcG_{f}^{1})^{\alpha} = 0\quad\hbox{\it for}\,\,\,
\alpha\in\bZ.
\leqno(1.6.2)
$$

\medskip\noindent
{\it Proof.}
It is well known that the action of the monodromy on the
$ 0 $-th Milnor cohomology
$ H^{0}(F_{x},\bC) $ is semisimple, and the multiplicity
of every eigenvalue is
$ 1 $.
(This can be reduced to the normal crossing case easily.)
So we get (1.6.2).
This implies (1.6.1) because
$ \bC\{t\}^{\alpha} = 0 $ for
$ \alpha\notin\bZ $.

\medskip
As a corollary, we get the following

\medskip\noindent
{\bf 1.7.~Proposition.}
{\it There is a canonical action of
$ \bC\lbb\rd_{t}^{-1}\rbb $ on
$ \cG_{f}^{i} $ for any
$ i $.
For
$ i = 1 $, it is compatible by {\rm (1.6.1)} with the
canonical action of
$ \bC\lbb\rd_{t}^{-1}\rbb $ on
$ \bC\{t\} $ which is defined by the integration without
constant term.
}

\medskip\noindent
{\it Proof.}
This follows from Proposition (1.4) if
$ i \ne 1 $, and we use also Lemma (1.6) for
$ i = 1 $.

\bigskip\bigskip
\centerline{\bf 2. Brieskorn modules}

\bigskip\noindent
{\bf 2.1.~Sheaf of Brieskorn modules.}
With the notation of (1.1), let
$ \cA_{f}^{\ssbull} $ be the complex whose
$ i $-th term is
$ \Ker(df\wdg : \Omega_{X}^{i} \to \Omega_{X}^{i+1}) $ and
whose differential is induced by
$ d $,
see [4], [10].
There is a natural inclusion
$$
\cA_{f}^{\ssbull} \to \cK_{f}^{\ssbull}.
\leqno(2.1.1)
$$
Let
$$
\cH_{f}^{i} = \cH^{i}\cA_{f}^{\ssbull}|_{X_{0}}.
$$
Note that
$ \cH^{i}\cA_{f}^{\ssbull} $ is supported on
$ X_{0} $ for
$ i \ne 0 $, and we have to take the restriction to
$ X_{0} $ in order to define the action of
$ \rd_{t}^{-1} $ on
$ \cH^{1}\cA_{f}^{\ssbull} $.
We will call
$ \cH_{f}^{i} $ the {\it sheaf of Brieskorn modules}
associated to
$ f $.
By (2.1.1) we have natural morphisms
$$
\cH_{f}^{i} \to \cG_{f}^{i}.
\leqno(2.1.2)
$$

Since the differential is
$ f^{-1}\cO_{S} $-linear,
$ \cH_{f}^{i} $ has a structure of a
$ \bC\{t\} $-module.
We have the action of
$ \rd_{t}^{-1} $ on
$ \cH_{f}^{i} $ by
$$
\rd_{t}^{-1}\omega = df\wdg \eta \quad
\hbox{with}\quad d\eta = \omega.
\leqno(2.1.3)
$$
This action is well defined by an argument similar to [4].
Indeed, for
$ i \ne 1 $, the ambiguity of
$ \eta $ is given by
$ d\eta' $ and
$ df\wdg d\eta' = -d(df\wdg\eta') $.
For
$ i = 1 $, the ambiguity of
$ \eta \in \cO_{X} $ is given by
$ \bC $, and we have a canonical choice of
$ \eta $ assuming that the restriction of
$ \eta $ to
$ X_{0} $ vanishes (this is allowed because
$ df\wdg d\eta = 0 $).

For
$ i = 1 $, we see that
$ \omega_{0} := df\otimes 1 $ in (1.1) belongs to
$ \cH_{f}^{1} $.
Here we can easily verify that
$ \cH_{f}^{1} \to \cG_{f}^{1} $ is injective, see (2.4) below.
We define the {\it sheaf of reduced Brieskorn modules}
$ \tcH_{f}^{i} $ by
$$
\tcH_{f}^{i} = \cH_{f}^{i}\,\,\,
\hbox{\rm if}\,\,\, i \ne 1,\,\,\,\hbox{\rm and}\,\,\,
\cH_{f}^{1}/\bC\{t\}\omega_{0} \,\,\,
\hbox{\rm if}\,\,\, i = 1.
\leqno(2.1.4)
$$

Let
$ \Omega_{X/S}^{\ssbull} $ denote the sheaf complex of
relative differential forms in the usual sense, i.e.
$ \Omega_{X/S}^{i} = \Omega_{X}^{i}/df\wdg\Omega_{X}^{i-1} $.
Then we have a canonical
$ f^{-1}\cO_{S} $-linear morphism
$$
df\wdg : \Omega_{X/S}^{\ssbull} \to \cA_{f}^{\ssbull}[1],
\leqno(2.1.5)
$$
which induces an isomorphism of complexes after the localization
by
$ f $.
Note that the image of
$ 1 $ by this morphism is
$ \omega_{0} $.

\medskip\noindent
{\bf 2.2.~Proposition.}
{\it Assume
$ X_{0} $ is a divisor with normal crossings.
Then
$ \cH_{f}^{i} $ and
$ \tcH_{f}^{i} $ are constructible sheaves of finite free
$ \bC\{t\} $-modules which are stable by the action of
$ t\rd_{t} $ and the eigenvalues of the residue of
$ t\rd_{t} $ are contained in
$ (-1,0] $,
see {\rm [5]}.
Furthermore, the canonical morphism
$ \cH_{f}^{i} \to \cG_{f}^{i} $ induces
isomorphisms compatible with the action of
$ \bC\{t\} $
$$
\cH_{f}^{i} = V^{>-1}\cG_{f}^{i},\quad
\tcH_{f}^{i} = V^{>-1}\tcG_{f}^{i},
\leqno(2.2.1)
$$
where
$ V $ is as in {\rm (1.3.3)}.
}

\medskip\noindent
{\it Proof.}
We first recall the proof of the corresponding assertion for
the relative logarithmic de Rham complex
$ \Omega_{X/S}^{\ssbull}(\log X_{0}) $,
see [21].
Let
$ (x_{1}, \dots, x_{n}) $ be a local coordinate system such
that
$ f = \prod_{i=1}^{r}x_{i}^{m_{i}} $.
Then
$$
df/f = \msum_{i=1}^{r}m_{i}dx_{i}/x_{i},
$$
and hence
$ (\Omega_{X,0}^{\ssbull}(\log X_{0}), df\wdg) $ is
acyclic, because it is identified with the Koszul complex
associated to the morphisms
$ m_{i} : \bC\{x\} \to \bC\{x\} $ for
$ 1\le i\le n $ where we put
$ m_{i} = 0 $ for
$ i > r $ (and
$ m_{i} \ne 0 $ for
$ i \le r) $.

Let
$ \rd_{i} = \rd/\rd x_{i} $,
$ \eta_{i} = dx_{i}/x_{i} $, and
$ \eta_{I} = \eta_{i_{1}}\wdg \cdots \wdg \eta_{i_{p}} $ for
$ I = \{i_{1}, \dots, i_{p}\} \subset \{2, \dots, n\} $.
Since
$ \eta_{i} \,(i > 1) $ and
$ df/f $ form a basis of
$ \Omega_{X,0}^{\ssbull}(\log X_{0}) $,
we see that the
$ \eta_{I} $ for
$ I \subset \{2, \dots, n\} $ form a basis of
$ \Omega_{X/S,0}^{p}(\log X_{0}) $ over
$ \bC\{x\} $.
Since
$ (\Omega_{X,0}^{\ssbull}(\log X_{0}), d) $ is
identified with the Koszul complex for the morphisms
$ x_{i}\rd_{i} \,(i\le r), \rd_{i} \,(i>r) $,
we can identify
$ (\Omega_{X/S,0}^{\ssbull}(\log X_{0}), d) $ with the
Koszul complex for
$$
x_{i}\rd_{i} - (m_{i}/m_{1})x_{1}\rd_{1}\,\,
(2\le i\le r), \,\,\, \rd_{i} \,\,(i>r).
$$
The last complex has a structure of double complex, and we may
assume
$ r = n $ by taking first the cohomology of the differential
defined by
$ \rd_{i} \,(i>r) $.

Let
$ e $ be the greatest common divisor of the
$ m_{i} \,(1\le i\le r) $,
and put
$ \mu_{i} = m_{i}/e $.
Then the morphisms
$ x_{i}\rd_{i} - (m_{i}/m_{1})x_{1}\rd_{1} $
preserve
$ \bC x^{\nu} $ for any
$ \nu = (\nu_{1}, \dots, \nu_{r}) \in \bN^{r} $,
and
$ H^{p}\Omega_{X/S,0}^{\ssbull}(\log X_{0}) $ is a free
$ \bC\{t\} $-module generated by
$$
x^{k\mu}\eta_{I}\,\,\,(0 \le k < e, \,\,|I| = p).
\leqno(2.2.2)
$$
For
$ \omega \in H^{p}\Omega_{X/S,0}^{\ssbull}(\log X_{0}) $,
put
$ \eta = t\rd_{t}\omega $.
It is the image of
$ \omega $ by the logarithmic Gauss-Manin connection, and
$ d\omega = f^{-1}df\wdg \eta $ by definition.
Let
$ \xi $ be a holomorphic vector field such that
$ \xi f = f $.
Let
$ L_{\xi} $ and
$ \iota_{\xi} $ denote respectively the Lie derivation and
the interior product.
Then
$$
L_{\xi}\omega = t\rd_{t}\omega \quad\hbox{\rm in}\quad
H^{p}\Omega_{X/S,0}^{\ssbull}(\log X_{0}),
\leqno(2.2.3)
$$
because
$ L_{\xi}\omega - d\iota_{\xi}\omega = \iota_{\xi}d\omega =
\iota_{\xi}(f^{-1}df\wdg \eta) =
\eta - f^{-1}df\wdg \iota_{\xi}\eta $.
(This holds for any holomorphic function
$ f $ having a vector field
$ \xi $ such that
$ \xi f = f $.)
So we get the assertion for the relative logarithmic complex.

Let
$ \tcA_{f}^{i} = \Ker(f^{-1}df\wdg : \Omega_{X}^{i}(\log X_{0})
\to \Omega_{X}^{i+1}(\log X_{0})) $.
By the acyclicity of the Koszul complex
$ (\Omega_{X,0}^{\ssbull}(\log X_{0}), f^{-1}df\wdg) $,
we have an isomorphism of complexes
$$
f^{-1}df\wdg : \Omega_{X/S,0}^{\ssbull}(\log X_{0}) \to
\tcA_{f,0}^{\ssbull}[1],
$$
and the cohomology of
$ \tcA_{f}^{\ssbull} $ gives the Deligne extension
such that the eigenvalues of the residue of the connection are
contained in
$ [-1,0) $.
Put
$ g = \prod_{i = 1}^{r}x_{i} $.
By the above calculation of the relative logarithmic complex
in (2.2.2--3), it is enough to show
$$
\cA_{f}^{i} = g \tcA_{f}^{i}.
$$
Here the inclusion
$ \cA_{f}^{i} \supset g \tcA_{f}^{i} $ is clear.
The opposite inclusion follows from the fact that a
meromorphic form
$ \omega $ is logarithmic if
$ g\omega $ is holomorphic and
$ f^{-1}df\wdg \omega $ is logarithmic.
(Indeed, we have
$ \omega = \sum_{I}a_{I}\eta_{I} $ where
$ a_{I} \in g^{-1}\bC\{x\} $ and
$ a_{I} $ does not have a pole along
$ x_{i} = 0 $ for
$ i \in I $.
If
$ a_{I} $ has a pole along
$ x_{j} = 0 $ with
$ j \notin I $, then
$ f^{-1}df\wdg a_{I}\eta_{I} $ is not logarithmic along
$ x_{j} = 0 $, and the highest order part of the pole of
$ x_{j}^{-1}dx_{j}\wdg a_{I}\eta_{I} $ along
$ x_{j} = 0 $ does not come from other
$ f^{-1}df\wdg a_{J}\eta_{J} $ with
$ J \subset I\cup\{j\} $.)

The compatibility with
$ \rd_{t}^{-1} $ is clear by the definition of the
differential of
$ \cK_{f}^{\ssbull} $.
This completes the proof of Proposition (2.2).

\medskip
The following gives a refinement of Hamm's result [10] and
is closely related to positivity (or negativity) theorems
in [11], [14].

\medskip\noindent
{\bf 2.3.~Theorem.}
{\it The
$ t $-torsion part of
$ \cH_{f}^{i} $ is annihilated by
$ t^{p} $ for
$ p \gg 0 $ locally on
$ X_{0} $,
and coincides with the kernel of the canonical morphism
$ \cH_{f}^{i} \to \cG_{f}^{i} $.
It coincides further with the
$ \rd_{t}^{-1} $-torsion part of
$ \cH_{f}^{i} $ for any
$ i $, and vanishes for
$ i = 1 $.
The stalks of the image of the morphism
$ \cH_{f}^{i} \to \cG_{f}^{i} $ are finite free modules over
$ \bC\{t\} $, which are contained in
$ V^{>-1}\cG_{f}^{i} $ and contain
$ V^{\alpha}\cG_{f}^{i} $ for
$ \alpha \gg 0 $.
}

\medskip\noindent
{\it Proof.}
Let
$ \pi : \tX \to X $ be an embedded resolution of
$ X_{0} $ as in the proof of (1.4).
Put
$ \tf = f\pi $.
By (2.2) together with a spectral sequence, the
$ \boR^{i}\pi_{*}\cA_{\tf}^{\ssbull} $ are
constructible sheaves of finite free
$ \bC\{t\} $-modules which are stable by the action of
$ t\rd_{t} $ and the eigenvalues of the residue of
$ t\rd_{t} $ are contained in
$ (-1,0] $.
Furthermore the kernel and the cokernel of the canonical
morphism
$$
\cH_{f}^{i} = \cH^{i}\cA_{f}^{\ssbull}|_{X_{0}}
\to \boR^{i}\pi_{*}\cA_{\tf}^{\ssbull}|_{X_{0}}
$$
are annihilated by
$ t^{j} $ for
$ j \gg 0 $ locally on
$ X_{0} $,
because the cohomology sheaves of the mapping cone of
$ \cA_{f}^{i} \to \boR\pi_{*}\cA_{\tf}^{i} $ for each
$ i $ are annihilated by
$ t^{j} $ for
$ j \gg 0 $ locally on
$ X_{0} $, using the coherence of the direct image sheaves.
So the
$ t $-torsion part of
$ \cH_{f}^{i} $ is annihilated by a high power of
$ t $ locally on
$ X_{0} $,
and the free part of
$ \cH_{f,x}^{i} $ is a finite free
$ \bC\{t\} $-module contained in
$ (\boR\pi_{*}\cA_{\tf}^{i})_{x} $.
(This is essentially the same argument as in [10],
see also [14].)
Furthermore,
$ (\boR\pi_{*}\cA_{\tf}^{i})_{x} $ is canonically
isomorphic to
$ V^{>-1}\cG_{f,x}^{i} $ by the theory of Deligne
extension [5] using a Milnor representative as in (1.5)
(or a spectral sequence as in (1.4.2)).

For
$ i \ne 1 $,
the action of
$ \rd_{t} $ is bijective on
$ \cG_{f,x}^{i} $ by Proposition 1.4, and
$ V^{>-1}\cG_{f,x}^{i} $ and
$ \cG_{f,x}^{i} $ are finite free modules over
$ \bC\lbb\rd_{t}^{-1}\rbb $ and
$ \bC\lbb\rd_{t}^{-1}\rbb[\rd_{t}] $ respectively,
see (1.3).
To show that the kernel of
$ \cH_{f}^{i} \to \cG_{f}^{i} $ coincides with the
$ \rd_{t}^{-1} $-torsion part, we have
$$
d(\msum_{j=0}^{r}\eta_{j}\otimes \rd^{j}) =
\omega\otimes 1,
$$
if and only if
$ d\eta_{0} = \omega $,
$ df\wdg\eta_{j} = d\eta_{j+1} \,(0 \le j < r) $ and
$ df\wdg\eta_{r} = 0 $.
Note that the first two equalities mean that
$ df\wdg\eta_{r} $ represents
$ \rd_{t}^{-r-1}\omega $ by (2.1.3).
These imply that the kernel of
$ \cH_{f}^{i} \to \cG_{f}^{i} $ is contained in the
$ \rd_{t}^{-1} $-torsion part.
Conversely, if the above first two equalities hold and
$ df\wdg\eta_{r} $ vanishes in
$ \cH_{f}^{i} $, then there exists
$ \eta_{r+1} $ such that
$ df\wdg\eta_{r} = d\eta_{r+1} $ and
$ df\wdg\eta_{r+1} = 0 $ by the definition of
$ \cH_{f}^{i} $.

For
$ i = 1 $, a similar argument shows that the canonical morphism
$ \cH_{f}^{1}\to\cG_{f}^{1} $ is injective and
the action of
$ \rd_{t}^{-1} $ on
$ \cH_{f}^{1} $ is injective, because
$ df\wdg : \cO_{X}\to\Omega_{X}^{1} $ is injective.
Then
$ \cH_{f}^{1} $ is
$ t $-torsion-free by the second assertion of this theorem.
This completes the proof of Theorem (2.3).

\medskip\noindent
{\bf 2.4.~Proposition.}
{\it The action of
$ \rd_{t}^{-1} $ on
$ \cH_{f}^{i} $ is naturally extended to that of
$ \bC\lbb\rd_{t}^{-1}\rbb $ on
$ \cH_{f}^{i} $.
This is compatible by the morphism
$ \cH_{f}^{i} \to \cG_{f}^{i} $ with the action of
$ \bC\lbb\rd_{t}^{-1}\rbb $ on
$ \cG_{f}^{i} $ in {\rm (1.7)}.
}

\medskip\noindent
{\it Proof.}
Assume first
$ i \ne 1 $ so that the action of
$ \rd_{t} $ on
$ \cG_{f}^{i} $ is bijective.
Then the canonical morphism
$ \cH_{f}^{i} \to \cG_{f}^{i} $ is compatible with the action of
$ \rd_{t}^{-1} $ by definition, and its image is stable by the
action of
$ \bC\lbb\rd_{t}^{-1}\rbb $ because it is a lattice stable by
$ \rd_{t}^{-1} $.
Using Lemma (2.5) below, there is a positive integer
$ p $ such that the torsion part of
$ \cH_{f}^{i} $ is annihilated by
$ \rd_{t}^{-p} $, and we have a canonical action of
$ \rd_{t}^{-p}\bC\lbb\rd_{t}^{-1}\rbb $ on
$ \cH_{f}^{i} $ compatible with the action of
$ \bC[\rd_{t}^{-1}] $ and also with the change of
$ p $.
Indeed, for
$ P \in \bC\lbb\rd_{t}^{-1}\rbb $ and
$ \omega \in \cH_{f}^{i} $, consider the image
$ \omega' $ of
$ \omega $ in
$ \cG_{f}^{i} $,
take a lifting of
$ P\omega' $ to
$ \cH_{f}^{i} $, and then multiply it by
$ \rd_{t}^{-p} $ so that the ambiguity of the lifting
is annihilated.
The compatibility follows from that of
$ \rd_{t}^{-1} $ with
$ \cH_{f}^{i} \to \cG_{f}^{i} $.

Now assume
$ i = 1 $.
With the notation of (1.3.3), we have
$$
\rd_{t}^{-1}\cG_{f}^{1} \subset V^{>0}\cG_{f}^{1}+
\tcG_{f}^{1},
\leqno(2.4.1)
$$
using the canonical splitting of (1.6.1).
On the other hand, we have
$$
\rd_{t}^{-1}\cH_{f}^{1} \subset V^{>0}\cG_{f}^{1},
\leqno(2.4.2)
$$
by the definition of
$ \rd_{t}^{-1}\omega $ in (2.1.3) which assumes the condition
$ \eta|_{X_{0}} = 0 $.
Indeed, (2.4.2) is proved, for example, by using the
period integral of
$ \eta $ along a horizontal family of topological
$ 0 $-cycles, see [14].
Note that the period integral of
$ \omega \in \cH_{f}^{i} $ is defined after dividing
$ \omega $ by
$ df $, i.e. by taking the inverse image by the morphism
(2.1.5).

Since the intersection of
$ V^{>0}\cG_{f}^{1} + \tcG_{f}^{1} $ with
$ \bC\omega_{0} $ vanishes, and the ambiguity of
$ \rd_{t}^{-1} $ is given by
$ \bC\omega_{0} $, the assertion follows.
Thus the proof of Proposition (2.4) is reduced to the following.

\medskip\noindent
{\bf 2.5.~Lemma.}
{\it
The
$ \rd_{t}^{-1} $-torsion part of
$ \cH_{f}^{i} $ is annihilated by
$ \rd_{t}^{-p} $ for
$ p \gg 0 $ locally on
$ X_{0} $.
}

\medskip\noindent
{\it Proof.}
By Theorem (2.3), this is reduced to the following:

\medskip\noindent
$ (A) \,\,\, $ Let
$ M $ be a torsion-free abelian group with an action of
$ t, \rd_{t}^{-1} $ satisfying the relation:
$ [t,\rd_{t}^{-1}] = \rd_{t}^{-2} $.
Assume
$ t^{p}M = 0 $ for a positive integer
$ p $.
Then
$ \rd_{t}^{-2p}M = 0 $.

\medskip
For this, we prove the following by increasing
induction on
$ q > 0 $:

\medskip\noindent
$ (A') \,\,\, $ If
$ m \in M $ satisfies
$ t^{q}m = 0 $ and
$ t^{p}M = 0 $, then
$ \rd_{t}^{-(p+q)}m = 0 $.

\medskip
Indeed,
$ t^{p+q-1}\rd_{t}^{-1}m $ is
a linear combination of
$ \rd^{-(p+q-i)}t^{i}m $ with
$ 0 \le i \le p + q - 1 $, and the coefficient of
$ \rd^{-(p+q)}m $ is
$ (p+q-1)! $ because
$ [t,\rd_{t}^{-j}] = j\rd_{t}^{-j-1} $ for
$ j > 0 $.
Then we can apply the inductive hypothesis to
$ t^{i}m $ with
$ q $ replaced by
$ q - i $ so that
$ \rd^{-(p+q-i)}t^{i}m = 0 $ for
$ i > 0 $.
Thus the assertion is proved.

\medskip\noindent
{\bf 2.6.~Remark.}
The above lemma follows also from a formula in [1] showing that
$ \rd_{t}^{-2p} $ is a linear combination of
$ \rd_{t}^{-j}t^{p}\rd_{t}^{-(p-j)} $.

\vfill\eject
\centerline{\bf 3. Finiteness property}

\bigskip\noindent
{\bf 3.1.~Triviality of
$ f $ along a stratification.}
Let
$ f $ be as in (1.1).
Let
$$
\Theta_{f} = \{\xi \in \Theta_{X} : \xi f = 0\},
$$
where
$ \Theta_{X} $ is the sheaf of holomorphic vector fields.
Let
$ \delta_{x} $ be the dimension of the image of
$ \Theta_{f,x} $ in the tangent space
$ T_{X,x} $.
Note that the restriction of
$ \cH_{f}^{i} $ to an integral curve of a nowhere vanishing vector
field
$ \xi\in\Theta_{f} $ is locally constant, see Lemma (3.2) below.
Let
$$
S_{i} = \{x \in X : \delta_{x} = i\}.
$$
We have
$ \dim S_{i} \ge i $ in general (by induction on
$ \dim X $ considering the integral curves of a vector field).
We say that
$ f $ is {\it locally trivial along a stratification} if

\medskip\noindent
$ (C) \,\,\,\, \dim S_{i} = i \,\, $ for any
$ i $.

\medskip\noindent
In this case, there is Whitney stratification
$ \{S'_{i}\} $ such that the tangent space of
$ S'_{i,x} $ is contained in the image of
$ \Theta_{f,x} $.
(This condition was considered in [17].)

Let
$ x \in X_{0} $,
and take a local coordinate system
$ (z_{1}, \dots, z_{n}) $ around
$ x $.
We have the distance function
$ \rho $ from
$ x $ defined by
$ (\sum_{i} |z_{i}|^{2})^{1/2} $ as usual.
If condition
$ (C) $ is satisfied, we see that the following condition is
satisfied (applying a standard argument to the above
stratification and
$ \rho) $:

\medskip\noindent
$ (C') \,\,\, $
For any
$ y \in X_{0} \setminus \{x\} $ near
$ x $,
there exists
$ \xi_{y} \in \Theta_{f,y} $ such that
$ \langle \xi_{y},d\rho\rangle_{y} \ne 0 $.

\medskip\noindent
This means that the integral curve of
$ \xi_{y} $ intersects transversely the sphere
$ S_{x,\varepsilon} $ with center
$ x $ and containing
$ y $ (i.e.
$ \varepsilon = \rho(y) $).

\medskip\noindent
{\bf 3.2.~Lemma.}
{\it Assume
$ f = \pi^{*}g $ with
$ \pi : X \to Y $ a smooth morphism of complex manifolds
{\rm (}i.e. the differential has the maximal rank{\rm )}.
Then
$ \cH_{f}^{i} = \pi^{-1}\cH_{g}^{i} $.
}

\medskip\noindent
{\it Proof.}
Since the assertion is local, we may assume
$ X = Y\times Z $ for a complex manifold
$ Z $.
The complexes
$ (\Omega_{X}^{\ssbull},df\wdg) $,
$ (\Omega_{X}^{\ssbull},d) $ and
$ (\cA_{f}^{\ssbull},d) $ have structures of double complexes, and
the differential for the second component of
$ (\Omega_{X}^{\ssbull},df\wdg) $ vanishes.
So
$ (\cA_{f}^{\ssbull},d) $ is the external product of
$ (\cA_{g}^{\ssbull},d) $ and
$ (\Omega_{Z}^{\ssbull},d) $.
Then, taking first the cohomology of the differential for the second
component and using the spectral sequence associated
to the double complex, the assertion follows.

\medskip

The following was essentially shown in [17] with
$ \cA_{f}^{\ssbull} $ replaced by
$ \Omega_{X/S}^{\ssbull} $.

\medskip\noindent
{\bf 3.3.~Theorem.}
{\it With the above notation, assume condition
$ (C) $ is satisfied.
Then
$ \cH_{f}^{i} $ is a constructible sheaf of finite modules over
$ \bC\{t\} $ and over
$ \bC\lbb\rd_{t}^{-1}\rbb $.
In particular, its torsion part is a constructible sheaf of finite
dimensional
$ \bC $-vector spaces.
}

\medskip\noindent
{\it Proof.}
By Lemma (3.2),
$ \cH_{f}^{i} $ is a constructible sheaf.
To show the finiteness, we may assume that
$ f : X \to S $ is a Milnor representative as in the proof
of (1.4) (in particular,
$ X $ is Stein).
Since the
$ \cA_{f}^{i} $ are coherent
$ \cO_{X} $-modules and the differential is
$ f^{-1}\cO_{S} $-linear, it is enough to show that
$ \boR^{i}f_{*}\cA_{f}^{\ssbull} $ is independent of
$ \varepsilon $ (shrinking
$ \delta $ if necessary) using the theory of
Kiehl-Verdier [13], see also [8].
So the assertion follows from the constructibility.

\medskip\noindent
{\bf 3.4.~Remark.}
We can show the finiteness of
$ \cH_{f,x}^{i} $ under the assumption
$ (C') $.

\medskip\noindent
{\bf 3.5.~Remark.}
Without condition
$ (C) $ or
$ (C') $, the torsion part of
$ \cH_{f,x}^{i} $ is not a finite dimensional
$ \bC $-vector space in general.
For example, consider
$ f = x^{5} + y^{5} + x^{3}y^{3}z $, or rather
$ f = x^{5}/5 + y^{5}/5 + x^{3}y^{3}z/3 $
replacing the coordinates.
Then
$$
\eqalign{
f_{x}
&= x^{4} + x^{2}y^{3}z = x^{2}(x^{2}+y^{3}z),\cr
f_{y}
&= y^{4} + x^{3}y^{2}z = y^{2}(y^{2}+x^{3}z),\cr
f_{z}
&= x^{3}y^{3}/3,\cr
}
$$
and
$ \Sing\,f = \{x = y = 0\} $.
Furthermore,
$ f $ is a topologically trivial deformation of holomorphic
functions of two variables, and the highest Milnor cohomology
$ H^{2}(F_{0},\bC) $ vanishes.
So it is enough to show that
$ \cH_{f,0}^{3} $ is not finite dimensional.

First we have to calculate
$ \cA_{f,0}^{2} \,(= \Ker\, df\wdg) $.
This is equivalent to the determination of the intersection of
the two ideals
$ (f_{x}, f_{y}) $ and
$ (f_{z}) $.
Since
$ f $ is weighted-homogeneous of weight
$ (1,1,-1) $, it is sufficient to consider only monomials
$ x^{i}y^{j}z^{k} $ with
$ i + j - k = c $ for some constant
$ c $.
So we may substitute
$ z = 1 $ for the calculation of the intersection of the ideals
(after taking the partial derivatives), and conclude that
$ \cA_{f,0}^{2} $ is generated over
$ \bC\{x,y,z\} $ by
$$
\eqalign{
&xy^{3} dy\wdg dz - 3(x^{2}+y^{3}z)dx\wdg dy,\cr
&x^{3}y\,dx\wdg dz + 3(y^{2}+x^{3}z)dx\wdg dy,\cr
&x^{2}y^{2}z\,dy\wdg dz + x^{3}dx\wdg dz +
3(y - xy^{2}z^{2})dx\wdg dy,\cr
&x^{2}y^{2}z\,dx\wdg dz + y^{3}dy\wdg dz -
3(x - x^{2}yz^{2})dx\wdg dy.\cr
}
$$
Then we see that the image of
$ d\cA_{f,0}^{2} $ is contained in the submodule
$ x\Omega_{X,0}^{3} + y\Omega_{X,0}^{3} $ which is
identified with the ideal
$ (x,y) $ of
$ \bC\{x,y,z\} $.
Since the action of
$ f $ on
$ \Omega_{X,0}^{3}/(x\Omega_{X,0}^{3} + y\Omega_{X,0}^{3}) $
is the multiplication by
$ 0 $, it implies that the torsion part of
$ \cH_{f,0}^{3} $ is infinite dimensional.
(A similar argument would apply also to
$ f = x^{5} + y^{5} + x^{2}y^{2}z $.)

\medskip\noindent
{\bf 3.6.~Remark.}
In [1], the following condition was considered for
$ i \ge 2 $:

\medskip\noindent
$ (P') \quad d(\Ker\, df\wedge) \cap \Im\,df\wedge =
\Im\,df\wedge d\quad\hbox{in}\quad\Omega_{X,x}^{i} $.

\medskip\noindent
We can easily show that this is a necessary and sufficient
condition for the torsion-freeness of
$ \cH_{f,x}^{i} $.
A stronger condition
$ (P) $ was also considered there in order to prove a
formula for the dimension of the Milnor cohomology.
In the case
$ \dim X = 2 $, it was proved in loc.~cit. that condition
$ (P) $ always holds, and hence also
$ (P') $ does:

\medskip\noindent
{\bf 3.7.~Proposition} [1].
{\it Assume
$ \dim X = 2 $.
Then condition
$ (P) $ always holds and hence the
$ \cH_{f}^{i} $ are torsion-free.
}

\medskip\noindent
{\bf 3.8.~Remark.}
We can show
$ (P') $ in the case
$ \dim X = 2 $, using a resolution of singularities as follows.
We may assume
$ i = 2 $.
Since
$ df\wedge $ is acyclic after localized by
$ f $,
any element of
$ \Ker\, df\wedge $ is uniquely written as
$ hdf/f $ with
$ h \in \cO_{X} $ such that
$ h|_{X_{0}} = 0 $.
Let
$ h_{x}, h_{y} $ denote the partial derivatives by local coordinates
$ x,y $.
Assume
$$
d(hdf/f) = dh\wdg df/f = (h_{x}f_{y}/f - h_{y}f_{x}/f)dx\wdg dy
\in \Im \,df\wedge.
\leqno(3.8.1)
$$
We will show that
$ h/f $ belongs to
$ \cO_{X} $ modifying
$ h $ by a linear combination of
$ f^{\alpha} $ where the
$ \alpha $ are rational numbers such that
$ f^{\alpha} $ is holomorphic (i.e. for any irreducible component
$ D_{i} $ of
$ f^{-1}(0) $ with multiplicity
$ m_{i} $, we have
$ m_{i}\alpha\in \bN $).
Note that
$ df^{\alpha}\wdg df/f = 0 $.

Let
$ \pi : \tX \to X $ be an embedded resolution of singularities of
$ f^{-1}(0) $.
Put
$ \tf = f\pi, \th = h\pi $ and
$ \tD = \tf^{-1}(0) $.
Note that the above conditions are compatible with the pull-back,
and it is sufficient to show that
$ \th/\tf $ is holomorphic on
$ \tX $, modifying
$ \th $ as above.

For a singular point
$ 0 $ of
$ \tD_{\rm red} $, let
$ x, y $ be local coordinates such that
$ \tf = x^{p}y^{q} $ around this point.
Let
$ V $ be the filtration on
$ \cO_{\tX,0}[\tf^{-1}] = \bC\{x,y\}[x^{-1}y^{-1}] $ such that
$ V^{\alpha} $ for
$ \alpha \in \bQ $ is generated over
$ \cO_{X} $ by
$ x^{i}y^{j} $ with
$ i/p \ge \alpha $,
$ j/q \ge \alpha $.
This can be defined also for a smooth point of
$ \tD_{\rm red} $ by considering only the condition
$ i/p \ge \alpha $ if
$ \tf = x^{p} $ locally.
So the filtration
$ V $ is globally well-defined on
$ \cO_{\tX}[\tf^{-1}] $, and
$ V^{1} $ is generated by
$ \tf $ over
$ \cO_{\tX} $.
Let
$ V^{>\alpha} = \mcup_{\beta >\alpha}V^{\beta} $.
The restriction of the graded pieces
$ \Gr_{V}^{\alpha} := V^{\alpha}/V^{>\alpha} $ to the smooth
points of
$ \tD $ is a line bundle if it does not vanish.

At a singular point of
$ \tD_{\rm red} $, let
$ g = xy(\th_{x}\tf_{y}/\tf - \th_{y}\tf_{x}/\tf) $.
Then
$ gdx\wdg dy/xy \in \Im\, \tf \wedge $ by (3.8.1) on
$ \tX $, and hence
$ g $ belongs to
$ V^{1} $ because
$ xy(\rd \tf) \subset V^{1} $ where
$ (\rd \tf) $ denotes the Jacobian ideal.
So we get
$$
qx\th_{x} - py\th_{y} \in V^{1}.
\leqno(3.8.2)
$$
Since the differential operator
$ qx{\rd\over\rd x} - py{\rd\over\rd y} $ preserves
$ \bC x^{i}y^{j} $ and
$ V^{\alpha} $,
we see that if
$ \th \in V^{\alpha} $ with
$ 0 < \alpha < 1 $ and (3.8.2) holds at every singular point of
$ \tD_{\rm red}, $ then there is a unique decomposition
$$
\th = \th' + c\tf^{\alpha}\quad \text{with}\quad \th' \in
V^{>\alpha}, \, c \in \bC,
\leqno(3.8.3)
$$
where
$ c = 0 $ unless
$ \tf^{\alpha} $ exists globally, i.e. unless
$ m\alpha \in \bN $ with
$ m $ the greatest common divisor of the multiplicities
$ m_{i} $ of
$ \tf $ along the irreducible components
$ \tD_{i} $ of
$ \tD $.
(Indeed, at singular points of
$ \tD_{\rm red} $, this follows easily from (3.8.2).
It implies that if the image of
$ \th $ in
$ \Gr_{V}^{\alpha} $ does not vanish generically on an irreducible
component
$ \tD_{i} $ of
$ \tD $ and hence
$ m_{i}\alpha $ is an integer,
then this holds also for the other components
$ \tD_{j} $ intersecting
$ \tD_{i} $ because otherwise
$ c = 0 $ at the intersection point.
Here we use analytic continuation on the smooth points of
$ \tD_{\rm red} $ together with the connectivity of
$ \pi^{-1}(0) $.)
So we may replace
$ \th $ with
$ \th' $,
and the assertion follows by increasing induction on
$ \alpha $.

\bigskip\bigskip
\centerline{\bf 4. Thom-Sebastiani type theorem}

\bigskip\noindent
{\bf 4.1.~External product.}
Let
$ f, g $ be nonconstant holomorphic functions on complex
manifolds
$ X $ and
$ Y $ respectively.
Define
$ h = f + g $ on
$ Z := X\times Y $.
Let
$ z = (x,y) \in X_{0}\times Y_{0} $,
and take
$ \omega \in \Omega_{X,x}^{i} $ such that
$ df\wdg \omega = d\omega = 0 $ with
$ i > 0 $.
Then we have the morphisms of complexes
$$
\cA_{g,y}^{\ssbull} \to \cA_{h,z}^{\ssbull}[i],\quad
\cK_{g,y}^{\ssbull} \to \cK_{h,z}^{\ssbull}[i],
$$
defined by
$ \eta \to \omega \wdg \eta $ and
$ \sum_{k} \eta_{k}\otimes \rd_{t}^{-k} \to
\sum_{k}\omega\wdg\eta_{k}\otimes \rd_{t}^{-k} $.
They induce morphisms
$$
\cH_{g,y}^{j} \to \cH_{h,z}^{i+j},\quad
\cG_{g,y}^{j} \to \cG_{h,z}^{i+j}.
\leqno(4.1.1)
$$
These are
$ \rd_{t}^{-1} $-linear if
$ j > 1 $.
For
$ j = 1 $,
we can show by increasing induction on
$ k $
$$
\omega \wdg g^{k}dg = 0\quad\text{in}\,\, \cH_{h}^{i+j}.
\leqno(4.1.2)
$$
Indeed, we have
$ \omega = d\omega' $ with
$ \omega' \in \Omega_{X,x}^{i-1} $,
and
$$
d(\omega'\wdg h^{k}dh) = \omega \wdg h^{k}dh =
\msum_{0\le i\le k}\hbox{${k\choose i}$}
f^{i}\omega \wdg g^{k-i}dg.
$$
So we can apply the inductive hypothesis to
$ f^{i}\omega \,\,(1\le i \le k) $ instead of
$ \omega $.

The assertion (4.1.2) means that the restriction of
(4.1.1) to
$ \bC[t]dg $ vanishes.
It implies that its restriction to
$ \bC\{t\}dg $ also vanishes, because the action of
$ \rd_{t}^{-1} $ on
$ \bC\{t\}/\bC[t] $ is bijective.
(Indeed, the last bijectivity implies that any
$ e \in \bC\{t\}/\bC[t] $ is infinitely
$ \rd_{t}^{-1} $-divisible, i.e. for any
$ r \in \bN $, there exists
$ e' \in \bC\{t\}/\bC[t] $ such that
$ \rd_{t}^{-r}e' = e $.
But this property does not hold for any
nonzero element of
$ \cH_{h}^{i+j} $ by (2.3) and (2.5).)
Thus (4.1.1) induces morphisms
$$
\tcH_{g,y}^{j} \to \tcH_{h,z}^{i+j},\quad
\tcG_{g,y}^{j} \to \tcG_{h,z}^{i+j},
\leqno(4.1.3)
$$
and these are
$ \bC\lbb\rd_{t}^{-1}\rbb $-linear and
$ \bC\lbb\rd_{t}^{-1}\rbb[\rd_{t}] $-linear respectively.
(Indeed, (2.3--5) implies that
$ \bigcap_{k\ge 0}\rd_{t}^{-k}\tcH_{h,z}^{j} = 0 $.)

We have an assertion similar to (4.1.3) with
$ \tcH_{g,y}^{j}, \tcG_{g,y}^{j} $ replaced by
$ \tcH_{f,x}^{i}, \tcG_{f,x}^{i} $,
if we take
$ \omega' \in \Omega_{Y,y}^{j} $ such that
$ dg\wdg \omega' = d\omega' = 0 $ with
$ j > 0 $.
So we get well-defined morphisms
$$
\eqalign{
\tcH_{f,x}^{i}\otimes_{\bC\lbb\rd_{t}^{-1}\rbb}
\tcH_{g,y}^{j}
&\to \tcH_{h,z}^{i+j},\cr
\tcG_{f,x}^{i}\otimes_{\bC\lbb\rd_{t}^{-1}\rbb[\rd_{t}]}
\tcG_{g,y}^{j}
&\to \tcG_{h,z}^{i+j},\cr
}
\leqno(4.1.4)
$$
such that the action of
$ t\otimes 1+1\otimes t $ on the left-hand side corresponds
to that of
$ t $ on the right-hand side.

\medskip\noindent
{\bf 4.2.~Theorem.}
{\it Let
$ n = \dim X $.
Assume
$ f $ has an isolated singularity at
$ x \in X_{0} $ so that
$ \tcH_{f,x}^{i} = \tcG_{f,x}^{i} = 0 $ for
$ i \ne n $.
Then, for
$ i = n $ and for any integers
$ j $, the morphisms {\rm (4.1.4)} induces isomorphisms
$$
\eqalign{
\tcH_{f,x}^{n}\otimes_{\bC\lbb\rd_{t}^{-1}\rbb}
\tcH_{g}^{j-n}
&\simto \tcH_{h}^{j}|_{\{x\}\times Y},\cr
\tcG_{f,x}^{n}\otimes_{\bC\lbb\rd_{t}^{-1}\rbb[\rd_{t}]}
\tcG_{g}^{j-n}
&\simto \tcG_{h}^{j}|_{\{x\}\times Y},\cr
}
\leqno(4.2.1)
$$
where
$ Y $ is identified with
$ \{x\}\times Y \subset Z $.}

\medskip\noindent
{\it Proof.}
Let
$ V $ be a finite dimensional
$ \bC $-vector subspace of
$ \Omega_{X,x}^{n} $ such that the morphism
$ V \to \Omega_{X,x}^{n}/df\wdg \Omega_{X,x}^{n-1} $ is
bijective.
Then
$$
\cH_{f,x}^{n} = \bC\lbb\rd_{t}^{-1}\rbb\otimes_{\bC}V,\quad
\cG_{f,x}^{n} = \bC\lbb\rd_{t}^{-1}\rbb[\rd_{t}]\otimes_{\bC}V,
$$
and there are injective morphisms of complexes
$$
\gamma : V\otimes_{\bC}\cA_{g,y}^{\ssbull}[-n] \to
\cA_{h,z}^{\ssbull},\quad
\gamma' : V\otimes_{\bC}\cK_{g,y}^{\ssbull}[-n] \to
\cK_{h,z}^{\ssbull}.
\leqno(4.2.1)
$$

We have to show that (4.2.1) induces isomorphisms of
reduced Brieskorn modules and of reduced Gauss-Manin systems.

Let
$ \ocA_{h,z}^{\ssbull} $,
$ \ocK_{h,z}^{\ssbull} $ denote the complexes defined
by their cokernels respectively.
Then the inclusion induces isomorphisms
$$
\beta :
\cH^{j}\ocA_{h,z}^{\ssbull} \simto
\cH^{j}\ocK_{h,z}^{\ssbull},
\leqno(4.2.2)
$$
i.e.
$ \ocK_{h,z}^{\ssbull}/\ocA_{h,z}^{\ssbull}
$ is acyclic.
Indeed,
$ \cK_{h,z}^{\ssbull} $ has the filtration
$ F $ defined by
$$
F_{p}\cK_{h,z}^{i} = \mopls_{k\le p+i}\Omega_{Z,z}^{i}
\otimes \rd_{t}^{k},
$$
and the graded pieces of its quotient filtration
$ F $ on
$ \ocK_{h,z}^{\ssbull}/\ocA_{h,z}^{\ssbull} $ are
$ \tau $-truncated complexes of
$ (\Omega_{Z,z}^{\ssbull}/V\otimes_{\bC}\Omega_{Y,y}^{\ssbull}
[-n], dh\wdg) $,
where
$ \tau $ denotes the canonical truncation which is defined by
using the kernel of the differential at a certain degree as
in [6].
Note that the Koszul complex
$ (\Omega_{Z,z}^{\ssbull}, dh\wdg) $ is the external
product of
$ (\Omega_{X,x}^{\ssbull}, df\wdg) $ and
$ (\Omega_{Y,y}^{\ssbull}, dg\wdg) $,
and is isomorphic to
$$
(\Omega_{Z,z}^{\ssbull}/V\otimes_{\bC}\Omega_{Y,y}^{\ssbull}
[-n], dh\wdg) \oplus
(V\otimes_{\bC}\Omega_{Y,y}^{\ssbull}[-n], dg\wdg),
$$
because
$$
(\Omega_{X,x}^{\ssbull}, df\wdg) =
(\Omega_{X,x}^{\ssbull}/V[-n], df\wdg) \oplus V[-n].
$$
Since
$ (\Omega_{X,x}^{\ssbull}/V[-n], df\wdg) $ is
acyclic, the above arguments imply the filtered acyclicity of
$ \ocK_{h,z}^{\ssbull}/\ocA_{h,z}^{\ssbull} $ (i.e.
the graded pieces are acyclic).

We now show that (4.2.1) induces isomorphisms
$$
\tgam_{j} : V\otimes_{\bC}\tcH_{g,y}^{j-n} \simto
\tcH_{h,z}^{j},\quad
\tgam'_{j} : V\otimes_{\bC}\tcG_{g,y}^{j-n} \simto
\tcG_{h,z}^{j}.
\leqno(4.2.3)
$$
By definition, the isomorphism
$ \beta $ in (4.2.2) is embedded into
the commutative diagram
$$
\matrix{
0 & \to & \Coker\,\cH^{j}\gamma & \to &
\cH^{j}\ocA_{h,z}^{\ssbull} & \to &
\Ker\,\cH^{j+1}\gamma & \to & 0\cr
&& \,\,\,\downarrow \scriptstyle{\beta'}
&& \,\,\downarrow \scriptstyle{\beta}
&& \,\,\,\downarrow \scriptstyle{\beta''} && \cr
0 & \to & \Coker\,\cH^{j}\gamma' & \to &
\cH^{j}\ocK_{h,z}^{\ssbull} & \to &
\Ker\,\cH^{j+1}\gamma' & \to & 0 \cr
}
$$
By (4.1.2--3) we may replace
$ \cH^{\ssbull}\gamma, \cH^{\ssbull}\gamma' $ in this
diagram with
$ \tgam_{\ssbull}, \tgam'_{\ssbull} $ in (4.2.3)
(i.e. the Brieskorn modules and the Gauss-Manin systems
are replaced with the reduced ones) by modifying
$ \cH^{j}\ocA_{h,z}^{\ssbull},
\cH^{j}\ocK_{h,z}^{\ssbull} $ appropriately,
but the middle vertical morphism
$ \tbet $ is still bijective,
where the vertical morphisms
$ \beta' $,
$ \beta $,
$ \beta'' $ for the new diagram are denoted respectively by
$ \tbet' $,
$ \tbet $,
$ \tbet'' $.
Note that (4.2.3) is equivalent to the vanishing
of the bottom row of the new diagram, because
$ \tbet $ is bijective.

Let
$ \cC' $ denote the category of
$ \bC\lbb\rd_{t}^{-1}\rbb $-modules which are annihilated
by sufficiently high powers of
$ \rd_{t}^{-1} $.
Considering the commutative diagram modulo
$ \cC' $, we may essentially neglect the torsion part of
the Brieskorn modules by (2.5).
Note that
$ \tgam'_{j} $ in (4.2.3) is obtained by the tensor of
$ \tgam_{j} $ with
$ \bC\lbb\rd_{t}^{-1}\rbb[\rd_{t}] $ over
$ \bC\lbb\rd_{t}^{-1}\rbb $, and this tensor is an exact functor,
and hence commutes with
$ \Ker $ and
$ \Coker $.
Thus the target of
$ \tbet' $ and
$ \tbet'' $
is obtained by the tensor of their source with
$ \bC\lbb\rd_{t}^{-1}\rbb[\rd_{t}] $
over
$ \bC\lbb\rd_{t}^{-1}\rbb $, and their cokernel is
obtained by the tensor with
$ \bC\lbb\rd_{t}^{-1}\rbb[\rd_{t}]/\bC\lbb\rd_{t}^{-1}\rbb $
over
$ \bC\lbb\rd_{t}^{-1}\rbb $ using the right exactness of tensor.
In particular, the vanishing of the target of
$ \tbet' $,
$ \tbet'' $ is reduced to that of their cokernels.
Since the cokernels are infinitely
$ \rd_{t}^{-1} $-divisible and does not belong to
$ \cC' $ if it does not vanish, the snake lemma implies
that each term of the bottom exact sequence of the
commutative diagram vanishes.
So the assertion follows.

\vfill\eject
\centerline{\bf References}

\bigskip

\item{[1]}
D.~Barlet,
Sur certaines singularit\'es non isol\'ees d'hypersurfaces I,
preprint, Feb. 2004.

\item{[2]}
A.~Beilinson, J.~Bernstein and P.~Deligne, Faisceaux Pervers,
Ast\'erisque, vol.~100, Soc.~Math.~France, Paris, 1982.

\item{[3]}
P.~Bonnet and A.~Dimca,
Relative differential forms and complex polynomials,
Bull. Sci. Math. 124 (2000), 557--571.

\item{[4]}
E.~Brieskorn, Die Monodromie der isolierten Singularit\"aten von
Hyperfl\"achen, Manu\-scripta Math., 2 (1970), 103--161.

\item{[5]}
P.~Deligne,
Equations Diff\'erentielles \`a Points Singuliers R\'eguliers,
Lect. Notes in Math. vol.~163, Springer, Berlin, 1970.

\item{[6]}
P.~Deligne, Th\'eorie de Hodge I, Actes Congr\`es Intern.
Math., Part 1 (1970), 425--430; II, Publ. Math. IHES, 40 (1971),
5--58; III, ibid. 44 (1974), 5--77.

\item{[7]}
P.~Deligne, Le formalisme des cycles \'evanescents, in
SGA7 XIII and XIV, Lect. Notes in Math. vol. 340, Springer,
Berlin, 1973, pp. 82--115 and 116--164.

\item{[8]}
A.~Douady, Le th\'eor\`eme des images directes de Grauert
(d'apr\`es Kiehl-Verdier), Ast\'erisque 16 (1974), 49--62.

\item{[9]}
H.A.~Hamm,
Ein Beispiel zur Berechnung der Picard-Lefschetz-Monodromie
f\"ur nicht\-isolierte Hyperfl\"achensingularit\"aten,
Math. Ann. 214 (1975), 221--234.

\item{[10]}
H.A.~Hamm, Habilitationsschrift.

\item{[11]}
M.~Kashiwara, $ B $-functions and holonomic systems,
Inv. Math. 38 (1976/77), 33--53.

\item{[12]}
M.~Kashiwara, Vanishing cycle sheaves and holonomic systems
of differential equations, Algebraic geometry (Tokyo/Kyoto,
1982), Lect. Notes in Math. 1016, Springer, Berlin,
1983, pp. 134--142.

\item{[13]}
R.~Kiehl and J.-L.~Verdier,
Ein einfacher Beweis des Koh\"arenzsatzes von Grauert,
Math. Ann. 195 (1971), 24--50.

\item{[14]}
B.~Malgrange, Int\'egrales asymptotiques et monodromie,
Ann. Sci. Ecole Norm. Sup. (4) 7 (1974), 405--430.

\item{[15]}
B.~Malgrange, Polyn\^ome de Bernstein-Sato et cohomologie
\'evanescente, Analysis and topology on singular spaces, II,
III (Luminy, 1981), Ast\'erisque 101--102 (1983), 243--267.

\item{[16]}
F.~Pham, Singularit\'es des Syst\`emes Diff\'erentiels de
Gauss-Manin, Progress in Math. vol. 2, Birkh\"auser, Basel,
1979.

\item{[17]}
M.~Saito, Gauss-Manin connection for non-isolated hypersurface
singularities (in Japanese),
master thesis, Feb. 1980, Univ. Tokyo.

\item{[18]}
M.~Saito, Hodge filtrations on Gauss-Manin systems, I,
J. Fac. Sci. Univ. Tokyo Sect. IA Math. 30 (1984), 489--498;
II, Proc. Japan Acad. Ser. A Math. Sci. 59 (1983), 37--40.

\item{[19]}
K.~Sakamoto, Milnor fiberings and their characteristic maps,
Manifolds--Tokyo 1973 (Proc. Internat. Conf., Tokyo, 1973),
Univ. Tokyo Press, Tokyo, 1975, pp. 145--150.

\item{[20]}
M.~Sebastiani and R.~Thom,
Un r\'esultat sur la monodromie, Inv. Math. 13 (1971), 90--96.

\item{[21]}
J.H.M.~Steenbrink, Limits of Hodge structures, Inv. Math. 31
(1975/76), no. 3, 229--257.

\bye